\documentclass{amsart}
\usepackage{amssymb}
\usepackage{enumerate}
\usepackage[active]{srcltx}

\newtheorem{theorem}{Theorem}
\newtheorem{claim}{Claim}[theorem]
\newtheorem{question}[theorem]{Question}

\theoremstyle{definition}

\newtheorem{definition}[theorem]{Definition}

\newcommand{\mc}[1]{\mathcal{#1}}
\newcommand{\mf}[1]{\mathfrak{#1}}

\newcommand{\setm}{\setminus}

\newcommand{\subs}{\subset}

\def\<{\left\langle}
\def\>{\right\rangle}
\def\br#1;#2;{\bigl[ {#1} \bigr]^ {#2} }

\newcommand{\alo}{{\aleph_{\omega}}}
\newcommand{\aloo}{{\aleph_{\omega+1}}}

\newcommand{\deltop}[1]{#1_{\delta}}

\author[L. Soukup]{Lajos Soukup}
\address
      { Alfr{\'e}d R{\'e}nyi Institute of Mathematics, 
        Hungarian Academy of Sciences}
\email{soukup@renyi.hu}
\urladdr{http://www.math-inst.hu/$\tilde{}$soukup}
\title[Noetherian type]
   {A note on   Noetherian type of spaces}
\date{\today}
\begin{document}
\begin{abstract}The {\em Noetherian type} of a space $X$, $Nt(X)$, is 
the least cardinal $\kappa$ such that $X$ has a base $\mc B$ such that 
$|\{B'\in \mc B: B\subs B'\}|<\kappa$ for each $B\in \mc B$.

Denote $\deltop{2^\alo}$ the space obtained from  $2^\alo$ by declaring the  $G_\delta$ sets to be open.
Milovich proved that if  
$\square_{\aleph_\omega}$ holds and $(\aleph_\omega)^\omega=\aleph_{\omega+1}$  
then $Nt(\deltop {2^{\aleph_\omega}})=\omega_1$.
Answering a question of Spadaro, we show that if $(\aleph_\omega)^\omega=\aleph_{\omega+1}$  
and  $(\aleph_{\omega+1}, \aleph_{\omega})\twoheadrightarrow (\aleph_1,\aleph_0)$
 then  $Nt(\deltop {2^{\aleph_\omega}})\ge\omega_2$.
\end{abstract}
\maketitle


\begin{definition}
 Let $X$ be a topological space. \\
(1) The {\em Noetherian type} of $X$, $Nt(X)$, is 
the least cardinal $\kappa$ such that $X$ has a base $\mc B$ such that 
$|\{B'\in \mc B: B\subs B'\}|<\kappa$ for each $B\in \mc B$.\\
(2)The {\em local Noetherian type} of
a point $x$ in  $X$, $\chi Nt(x,X)$, is 
the least cardinal $\kappa$ such that $x$ has a local base $\mc B$  in $X$ such that 
$|\{B'\in \mc B: B\subs B'\}|<\kappa$ for each $B\in \mc B$.
The  {\em local Noetherian type} of
a space $X$, $\chi Nt(X)$, is defined by the formula
$\chi Nt(X)=\sup\{\chi Nt(x,X):x\in X\}$. 
\end{definition}

For a topological space $X$, let $\deltop X$ denote  
the space obtained by declaring the $G_\delta$-sets to be open.

\begin{theorem}[Spadaro]\label{spadaro} (GCH) Let X be a compact space such that $Nt(X)$ has
uncountable cofinality. Then $Nt(\deltop X) \le 2^{N t(X)}$.  
\end{theorem}

He asked what happens if $Nt(X)$ has countable cofinality or if we drop GCH.
Since Milovich, \cite{36}, proved that if  $X$  is  compact dyadic homogeneous space then $Nt(X)=\omega$,
 Theorem \ref{spadaro} does not  apply  to the Cantor cubes.
However Milovich also proved that $Nt (\deltop{2^\kappa})=\omega_1$ for $\kappa<\alo$.
So the  simplest unsettled case remained $\kappa=\alo$. 

Spadaro formulated the following question:

\begin{question}
Is $Nt(\deltop{2^{\alo}})=\mf c^+$? 
\end{question}

\begin{theorem}[Milovich]
 If  $\square_{\aleph_\omega}$ holds and $(\aleph_\omega)^\omega=\aleph_{\omega+1}$  
then $Nt(\deltop {2^{\aleph_\omega}})=\omega_1$.
\end{theorem}

We prove the following result:
\begin{theorem}\label{tm:chang}  If $(\aleph_\omega)^\omega=\aleph_{\omega+1}$  
and 
$(\aleph_{\omega+1}, \aleph_{\omega})\twoheadrightarrow (\aleph_1,\aleph_0)$
 then  $Nt(\deltop {2^{\aleph_\omega}})\ge\omega_2$.
\end{theorem}

Recall that 
$(\kappa,\lambda)\twoheadrightarrow(\mu,\nu)$ is the following statement: \\
{\em For any structure ${\mc A}=(A,U,\dots)$ of countable signature with 
$|A|=\kappa$, $U\subseteq A$ and $|{U}|=\lambda$, there is an 
elementary substructure ${\mc A}'=(A',U',\ldots)$ of $\mc A$  such that  
$| {A'}|=\mu$ and $| {U'}|=\nu$. 
}

In \cite{levinski}, a model of ZFC $+$ GCH $+$ Chang's Conjecture for 
$\aleph_\omega$, i.e.\ 
$(\aleph_{\omega+1},\aleph_\omega)\twoheadrightarrow(\aleph_1,\aleph_0)$, was
constructed starting from some quite strong large cardinal assumption.

Theorem \ref{tm:chang} will follow easily from Theorem \ref{tm:chang2} below.

\begin{theorem}\label{tm:chang2}  Assume that  $(\aleph_\omega)^\omega=\aleph_{\omega+1}$  
and  $(\aleph_{\omega+1}, \aleph_{\omega})\twoheadrightarrow (\aleph_1,\aleph_0)$.
If $\mc D
$ is cofinal in $\<\br \alo;\omega;,\subseteq\>$, then 
there is  $A\in \br \alo;\omega;$ such that  $|\mc D\cap \mc P(A)|> \omega$.
\end{theorem}

\begin{proof}[Proof of Theorem \ref{tm:chang2}]
Assume on  the contrary that $\mc D$ is cofinal, but 
$|\mc D\cap \mc P(A)|\le \omega$ for each $A\in \br \alo;\omega;$.

Let us fix an enumeration $\{a_\alpha:\alpha<\aloo\}$ of 
$[\alo]^{\aleph_0}$ without repetition.
Consider 
\begin{enumerate}[(1)]
\item the relation $E=\{\<\alpha,\beta\>:\alpha\in\alo,\beta\in\aloo,\alpha\in a_\beta\}$,
\item a function ${f}:\aloo\to {\aloo}$  such that  
$a_\alpha\subset a_{f(\alpha)}\in \mc D$ for each $\alpha<\aloo$,
\item a function ${g}:{\aloo\times\aloo}\to{\alo}$ such that, 
$g(\alpha,\cdot)\restriction\alpha$ is an injective mapping 
from $\alpha$ to $\alo$    for each 
$\alpha\in\aloo$,
\item a function ${h}:{\aloo\times\aloo\times\aloo}\to {\omega}$  such that  
for each $\alpha$, $\beta\in\aleph_{\omega+1}$, 
the function $h(\alpha,\beta,\cdot)$ is an bijection  between $a_\alpha\setminus a_\beta$ and 
an initial segment  of $\omega$,
\end{enumerate}
and consider the structure:
\begin{equation}\notag
\mc A=(\aloo, \alo,\leq, E,f,g,h),
\end{equation}
where $\le$ is the natural ordering on $\aloo$.

By (4) and since $\omega$ is definable in $\mc A$, we 
can express ``{\em $a_\alpha\setm a_\beta$ is 
infinite}'' as a formula $\varphi(\alpha,\beta)$ in the language of $\mc A$. 
Now, applying the Chang's conjecture for $\mc A$ with 
$A=\aleph_{\omega+1}$ and $U=\aleph_\omega$, we obtain 
elementary substructure $\mc A'=(A',U',\leq',E',f',g',h')$ of $\mc A$ 
 such that  $| {A'}|=\aleph_1$ and $| {U'}|=\aleph_0$. 
\begin{claim}\label{claim1}
$otp(A')=\omega_1$.
\end{claim}
\begin{proof}[Proof of the Claim]
By (3) and elementarity of $\mc A'$, every initial segment of 
$A'$ can be mapped into $U'$ injectively and hence countable. Since 
$| {A'}|=\aleph_1$, it follows that $otp(A')=\omega_1$.

\end{proof}

\begin{claim}\label{claim2}
For any $\alpha<\aloo$, there is $\gamma<\aloo$  such that  
$a_\gamma\setm a_\beta$ is infinite for every $\beta<\alpha$.
\end{claim}

\begin{proof}[Proof of the Claim]

Since $| {\alpha}|\leq\alo$, we can find a partition 
$\{I_n:{n\in\omega}\}$ of $\alpha$  such that  $| {I_n}|<\alo$ for every 
$n<\omega$. For $n<\omega$, let 
\begin{equation}\notag
\eta_n=\min\Bigl(\alo\setminus\bigl(\{\eta_i:i<n\}\cup\bigcup\{{a_\xi}:
{\xi\in\bigcup_{m\leq n}I_m}\}\bigr)\Bigr),
\end{equation}
and let $\gamma<\aloo$ be  such that  
$a_\gamma=\{{\eta_n}:{n\in\omega}\}$. For any $\beta<\alpha$, if $\beta\in I_{m}$ for some 
$m<\omega$, then we have 
$a_\beta\cap a_\gamma\subseteq\{{\eta_n}:{n<m}\}$. Thus this $\gamma$ 
is as desired.
\end{proof}

\begin{claim}\label{claim3}
For any countable $I\subseteq A'$, there is $\gamma\in A'$  such that  
$a_\gamma\setm a_\beta$ is infinite for every $\beta\in I$.
\end{claim}
\begin{proof}[Proof of the Claim]
By Claim \ref{claim1}, there is $\alpha\in A'$  such that  $I\subseteq\alpha$. By 
elementarity of $\mc A'$, the formula with the parameter $\alpha$ 
expressing the statement of Claim \ref{claim2} for this $\alpha$ holds in 
$\mc A'$. Hence there is some $\gamma\in A'$  such that  
$a_\gamma\setm a_\beta$ is infinite for every $\beta\in A'\cap\alpha$. 
\end{proof}
After this preparation we are ready to proof Theorem \ref{tm:chang2}. Let 
\begin{equation}\notag
 I=\{{f(\xi):\xi \in A'}\}. 
\end{equation}
Since $A'$ is $f$-closed, we have $I\subs U'$.
If $\zeta\in A'$ then $a_\zeta\subs U$, so 
$\{a_\zeta:\zeta\in I\}=\{a_{f(\xi)}:\xi\in A'\}\subs \mc P(U')$ for $\xi \in A'$.
So $|I|\le \omega$ because $a_{f(\xi)}\in \mc D$ and  $|\mc D\cap \mc P(U')|\le \omega$. Hence, by 
\ref{claim3}, there is $\gamma\in A'$  such that  $a_\gamma\setm a_{\beta'}$ is 
infinite for every $\beta'\in I$. 
Consider the ordinal $\beta=f(\gamma)$. 
Then $a_\gamma\setm a_\beta$ should be infinite because $\beta=f(\gamma)\in I$.  
On the other hand,  $a_\beta=a_{f(\gamma)}\supset a_\gamma$, 
so $a_\gamma\setminus a_\beta=\emptyset$.
This contradiction proves the theorem.
\end{proof}

\begin{proof}[Proof of Theorem \ref{tm:chang}]
We show that $\chi Nt(x,\deltop {2^{\aleph_\omega}})=\omega_2$
for each $x\in 2^{\aleph_\omega}$.

Assume to the contrary that  $x$ has a local base $\mc B$  in $X$ such that 
\begin{equation}\label{b}
\text{$|\{B'\in \mc B: B\subs B'\}|\le \omega$ for each $B\in \mc B$}. 
\end{equation}

We can assume that $B=[x\restriction d_B]$
for some $d_B\in \br \alo;\omega;$.
Consider the family $\mc D=\{d_B:B\in \mc B\}
$.
Then 
\begin{enumerate}[(a)]
\item $\mc D$ is cofinal in $\<\br \alo;\omega;,\subseteq\>$, 
\item $|\mc D\cap \mc P(A)|\le \omega$ for each $A\in \br \alo;\omega;$.
\end{enumerate}
Observe that $[x\restriction D]\subset [x\restriction A]$
iff $D\supset A$. So
(a) holds because $\mc B$ is a local base of $x$, and  \ref{b} implies  (b).
However, this contradicts  Theorem \ref{tm:chang2}.
\end{proof}

\end{document}